\documentclass[11pt, letterpaper]{article}

\usepackage{amsmath}
\usepackage{amssymb}
\usepackage{amsfonts}
\usepackage{bbm}

\newtheorem{theorem}{Theorem}
\newtheorem{proposition}{Proposition}

\newtheorem{definition}{Definition}

\newtheorem{claim}{Claim}

\newcommand{\PP}{\mathrm{P}}
\newcommand{\NP}{\mathrm{NP}}
\newcommand{\ch}{\mathrm{ch}}

\newcommand{\prob}{\mathbb{P}}

\newcommand{\ind}{\mathbbm{1}}

\newcommand{\bfp}{\mathbf{p}}
\newcommand{\ovp}{\bar{p}}

\newcommand{\lcal}{\mathcal{L}}

\newcommand{\eps}{\varepsilon}

\newcommand{\wtprob}{\widetilde{\prob}}
\newcommand{\tf}{\tilde f}

\newenvironment{proof}{\noindent{\textbf{Proof.}}} {$\blacksquare$\vskip \belowdisplayskip}

\begin{document}

\title{A Short Proof that Phylogenetic Tree
Reconstruction\\ by Maximum Likelihood is Hard}

\author{{\bf S\'ebastien Roch}\\
Department of Statistics\\
University of California, Berkeley}

\maketitle

\begin{abstract}
Maximum likelihood is one of the most widely used techniques to infer
evolutionary histories. Although it is thought to be 
intractable, a proof of its hardness has been 
lacking.
Here, we give a short proof that computing the maximum likelihood
tree is NP-hard
by exploiting a connection between likelihood
and parsimony observed by Tuffley and Steel.
\end{abstract}

\section{Introduction}

In a series of seminal works, Edwards and Cavalli-Sforza~\cite{EC64},
Neyman~\cite{Ne71}, and Felsenstein~\cite{Fe81} applied
the maximum likelihood methodology to the problem of inferring phylogenies
from molecular sequences. 
Since, the many variants of this approach have gained
increasing popularity in the systematics literature. This is due in part 
to the flexibility of the technique 
in accomodating a variety of models of evolution as well as
to good practical performance. Nevertheless, the approach is not without a flaw: 
it has been observed to be highly demanding computationally.  
Remarkably, the computational complexity status of the problem
has remained elusive.
Partial progress was made recently in~\cite{AC04} where a variant of the
problem, known as Ancestral Maximum Likelihood, was shown to be $\NP$-hard.
Here, we resolve the issue by proving that computing the maximum likelihood
tree is $\NP$-hard. Moreover, we show that the log-likelihood
is $\NP$-hard to approximate within
a constant ratio. Our proof---which is mostly elementary---combines a connection between likelihood
and parsimony observed by Tuffley and Steel~\cite{TS97} with a result on
the hardness of approximating parsimony obtained by Wareham~\cite{Wa93}.

General references on 
inferring phylogenies are~\cite{Fe03,SS03}. 
For background on $\NP$-completeness and
hardness of approximation, refer to~\cite{GJ76, AC99}.
Many other popular phylogenetic techniques have been shown
to be $\NP$-hard, including parsimony~\cite{FG82,DJS86},
compatibility~\cite{DS86}, 
and distance-based methods~\cite{AB+99}.

\textsc{Remark.} While writing this paper, 
Mike Steel brought to our attention
that Benny Chor and Tamir Tuller have recently given an independent proof of
this result which now appears in the proceedings of RECOMB 2005~\cite{CT05}. 
Similarly to the proof presented here,
the Chor-Tuller paper uses results from~\cite{TS97} and
the hardness of approximating vertex cover
(from which follows the hardness of approximating 
parsimony), 
but their argument proceeds from a sequence
of rather involved constructions.
Our reduction has the advantage
of being short and elementary. It also sheds some more light on the
interesting connection between likelihood and parsimony.

\section{Definitions and Results}

The use of maximum likelihood requires the choice of a statistical
model of evolution. Here, we consider the simple binary symmetric model
generally known as the Cavender-Farris model~\cite{Ca78, Fa73}.
We are given a tree $T$ on $n$ leaves and probabilities of transition
on edges $\bfp = \{p_e\}_{e\in E(T)} \in [0,1/2]^{E_T}$, where $E(T)$ is the set
of edges of $T$ and $E_T \equiv |E(T)|$ is the cardinality of $E(T)$. 
(All trees considered here have no internal vertex
of degree $2$.)
A realization of the model is obtained as follows: choose any vertex
as a root; pick a state for the root uniformly at random in $\{0,1\}$; moving
away from the root, each edge $e$ flips the state of its ancestor with probability
$p_e$. Let $[n]$ denote the set of leaves. 
A character $\chi$ assigns to each leaf a state
in $\{0,1\}$. An extension of $\chi$ is an assignment
of states in $\{0,1\}$ to all vertices of $T$ which is equal to $\chi$ on the leaves. 
The set of all extensions of $\chi$ is denoted $H(\chi)$.

In the Cavender-Farris model, the (modified) log-likelihood of $\chi$ on $(T,\bfp)$ 
is
\begin{equation}
\widetilde\lcal(\chi;T, \bfp) \equiv -\ln 2 \prob[\chi\,|\,T,\bfp]
= -\ln\left(\sum_{\hat\chi \in H(\chi)} \prod_{e=(u,v)\in E(T)}
p_e^{\ind\{\hat\chi(u)\neq\hat\chi(v)\}}(1-p_e)^{\ind\{\hat\chi(u)=\hat\chi(v)\}} \right),
\end{equation}
where $\ind\{A\}$ is $1$ if $A$ occurs, and $0$ otherwise.
The data consists of a set
of characters $X=\{\chi_i\}_{i=1}^{k}$.
Assuming the characters are independent and identically distributed,
the log-likelihood of the data is the sum of the log-likelihood of all
characters, viz.
$\widetilde\lcal(X;T, \bfp) = \sum_{i=1}^k \widetilde\lcal(\chi_i;T, \bfp)$.
The maximum likelihood (ML) problem consists in
computing $(T^*,\bfp^*)$ minimizing $\widetilde\lcal(X;T, \bfp)$ over
all trees and transition probability vectors. Note that this is equivalent
to minimizing the standard likelihood, without the factor of $2$.

Contrary to ML, maximum parsimony (MP) is not based on a
statistical model.
Denote by $\ch(\hat\chi)$ the number
of flips in $\hat\chi$, i.e.
$\ch(\hat\chi) = |\{(u,v) \in E(T)\,:\, \hat\chi(u)\neq\hat\chi(v)\}|$.
Let $l(\chi,T)$ be the smallest number of flips in any extension
of $\chi$ on $T$, i.e.
$l(\chi,T) = \min_{\hat\chi\in H(\chi)} \ch(\hat\chi)$.
The parsimony score of $T$ is then
$l(X,T) = \sum_{i=1}^k l(\chi_i, T)$.
The problem MP consists in finding
the tree $T^{**}$ minimizing $l(X, T)$ over all
trees. 

A useful connection between
ML and MP
was noted by Tuffley and Steel in~\cite{TS97}:
if one adds sufficiently many constant sites (i.e.
$\chi(i) = \alpha,\ \forall i\in [n]$ for some $\alpha\in\{0,1\}$) to the data and
applies the maximum likelihood technique, then
one necessarily chooses the most parsimonious tree. 
This could serve as the basis for
a reduction, except that the Tuffley-Steel bounds require
an exponential number of constant sites. Our
contribution is to show that a polynomial number of sites
imposes a weaker relationship between likelihood and parsimony,
but that this is sufficient for the following reason.
Parsimony is in fact \emph{hard to approximate},
that is, even the seemingly easier task of obtaining a solution
close to optimal is hard. This result is due to Wareham~\cite{Wa93}.

We prove the following theorem. We first define the notion of approximation algorithm.
\begin{definition}
Let $\Pi$ be an optimization problem (minimization). Let $I$ denote
an instance of $\Pi$ and $\mathrm{OPT}(I)$, the optimal value
of a solution to $I$. For $c>0$, a $(1+c)$-approximation algorithm
for $\Pi$ is a polynomial-time algorithm that is guaranteed to 
return, for any instance $I$, a solution with objective value
$m$ satisfying
$m \leq (1+c)\mathrm{OPT}(I)$.
\end{definition}
\begin{theorem}\label{theo:1}
There exists a $c>0$ sufficiently small so that there is no 
$(1+c)$-approximation algorithm for 
ML unless $\PP=\NP$. In particular, ML is $\NP$-hard.
(The approximation claim relates to the modified log-likelihood
$\widetilde\lcal$.)
\end{theorem}

\section{Proof}

In this section, we prove our main result.
The proof follows easily from the following propositions.
The first proposition borrows heavily
from~\cite{TS97} although we need somewhat tighter estimates.
The second proposition follows directly from the work of~\cite{Wa93,CT99}.
\begin{proposition}\label{prop:1}
Let $c' > c > 0$ be constants.
If there is a $(1+c)$-approximation algorithm
for ML then there is a $(1+c')$-approximation
algorithm for MP.
\end{proposition}
\begin{proposition}[\cite{Wa93,CT99}]\label{prop:2}
There exists a $c'>0$ sufficiently small so that there is no 
$(1+c')$-approximation algorithm for 
MP unless $\PP=\NP$.
\end{proposition}

As in~\cite{TS97}, the reduction from 
MP to ML
consists in adjoining
a large number of constant sites to the data.
Let $\eps> 0$ be a small constant 
and $M = \max\{2n,k\}$. Fix $N_c = M^{1/\eps}$. 
Denote by $X_0 = \{\chi_i\}_{i=1}^{k+N_c}$ the set $X$ augmented
with $N_c$ all-$0$ characters. 
For all $\chi$, let $N_\chi$ be the number of characters
equal to $\chi$ in $X$. To avoid the factors of $2$
from the probability of the root state, we define
$\wtprob[X_0\,|\,\bfp,T] = 2^{k+N_c} \prob[X_0\,|\,\bfp,T]$ for $\bfp\in[0,1/2]^{E_T}$.
Also, let $\tf_\chi = 2 \prob[\chi\,|\,\bfp,T]$.
Let $0$ be the all-$0$ character and $1$, the
all-$1$ character.
We make three claims, from which Proposition~\ref{prop:1}
follows. 
\begin{claim}\label{claim:1}
Let $\eps> 0$ and $N_c = M^{1/\eps}$ for $M = \max\{2n,k\}$.
Let 
$
p_e = q = \frac{l(X,T)}{E_T(k+N_c)}$,
for all $e\in E(T)$.
Then
$-\frac{\ln \wtprob[X_0\,|\,\bfp,T]}{\ln (k+N_c)} \leq (1+2\eps)l(X, T)$,
for $M$ large enough.
\end{claim}
\begin{proof}
Note that, by a calculation identical to~\cite[Lemma 5]{TS97},
\begin{equation*}
\ln \tf_\chi 
= \ln\left(\sum_{\hat\chi \in H(\chi)} q^{\ch(\hat\chi)} (1-q)^{E_T - \ch(\hat\chi)}\right)
\geq \ln\left( q^{l(\chi,T)} (1-q)^{E_T}\right)
\geq l(\chi,T)\ln q - E_T \left(q + 2 q^2\right),   
\end{equation*}
where we have used a 
standard Taylor expansion (note that we have $q\leq 1/2$ by definition).
This bound applies in particular to the case $\chi = 0$.
Then, as in~\cite[Lemma 5]{TS97} again,
\begin{eqnarray*}
-\frac{\ln \wtprob[X_0\,|\,\bfp,T]}{\ln (k+N_c)}
&=& -\frac{1}{\ln (k+N_c)}\ln\left(\tf_0^{N_0 + N_c} \prod_{\chi\neq 0}\tf_\chi^{N_\chi}\right)\\
&\leq& \frac{1}{\ln (k+N_c)}\left((k + N_c)E_T(q + 2q^2) - l(X,T)\ln q\right)\\
&=& l(X,T)\left(1 + \frac{\ln E_T - \ln l(X,T)}{\ln (k+N_c)} + \frac{1}{\ln (k+N_c)}\left(1 + 2 \frac{l(X,T)}{E_T(k+N_c)}\right)\right)\\
&\leq& l(X,T)\left(1 + \frac{\ln M}{\ln M^{1/\eps}} + \frac{1}{\ln M^{1/\eps}}\left(1 + 2 \frac{M}{M^{1/\eps}}\right)\right)\\
&\leq& (1 + 2\eps)l(X,T),
\end{eqnarray*}
for $M$ large enough. (Note that $\wtprob[X_0\,|\,\bfp,T] \leq 1$ because
$\prob[\chi\,|\,\bfp,T]\leq 1/2$ by symmetry.)
\end{proof}

\begin{claim}\label{claim:2}
For all $\bfp\in[0,1/2]^{E_T}$ such that
$-\frac{\ln \wtprob[X_0\,|\,\bfp,T]}{\ln (k+N_c)} \leq l(X, T)$,
one has $p_e \leq \ovp$, $\forall e\in E(T)$,
with $\ovp \equiv \frac{l(X,T)\ln (k+N_c)}{N_c}$.
\end{claim}
\begin{proof}
Assume edge $e$ is such that $p_e > \ovp$. Take
any two leaves $u,v$ joined by a path going through
$e$. As observed in~\cite[Formula (11)]{TS97}, the probability that
$\chi(u)\neq\chi(v)$ is at least $p_e$. In particular,
the probability that a character is constant is less
than $1 - p_e$ and
$-\ln \tf_0 \geq -\ln(1 - p_e) \geq p_e$ (by the $0-1$ symmetry).
Therefore,
\begin{eqnarray*}
-\frac{\ln \wtprob[X_0\,|\,\bfp,T]}{\ln (k+N_c)} 
= -\frac{1}{\ln (k+N_c)}\ln\left(\tf_0^{N_0 + N_c} \prod_{\chi\neq 0}\tf_\chi^{N_\chi}\right)
\geq -\frac{1}{\ln (k+N_c)}\left(\ln \tf_0^{N_c}\right)
> l(X,T),
\end{eqnarray*}
by $p_e > \ovp$ and $\tf_\chi \leq 1$ (by the $0-1$ symmetry). This contradicts the assumption.
\end{proof}

\begin{claim}\label{claim:3}
Let $\eps> 0$ and $N_c = M^{1/\eps}$ for $M = \max\{2n,k\}$.
For all $\bfp\in[0,1/2]^{E_T}$,  we have
$-\frac{\ln \wtprob[X_0\,|\,\bfp,T]}{\ln (k+N_c)} \geq (1-5\eps) l(X, T)$ for $M$ large enough.
\end{claim}
\begin{proof}
For this proof, we need a better estimate than~\cite[Lemma 6]{TS97}. 
From Claim~\ref{claim:2}, the result holds
whenever $\max_e p_e > \ovp$. Therefore,
we can assume that
for all $e\in E(T)$, $p_e \leq \ovp$. Then,
\begin{equation*}
\tf_\chi
\leq \sum_{\hat\chi \in H(\chi)} \ovp^{\ch(\hat\chi)}
\leq \sum_{\alpha = 0}^{E_T - l(\chi,T)} \binom{E_T}{\alpha + l(\chi,T)} \ovp^{\alpha + l(\chi,T)}
\leq \sum_{\alpha = 0}^{E_T - l(\chi,T)} (E_T\ovp)^{\alpha + l(\chi,T)}
\leq E_T (E_T\ovp)^{l(\chi,T)},
\end{equation*}
when $M$ is large enough so that $\ovp < 1/E_T$. For constant sites,
we use the bound $\tf_0, \tf_1\leq 1$ (by the $0-1$ symmetry). Therefore,
\begin{eqnarray*}
-\frac{\ln \wtprob[X_0\,|\,\bfp,T]}{\ln (k+N_c)}
&=& -\frac{1}{\ln (k+N_c)}\ln\left(\tf_0^{N_0 + N_c} \prod_{\chi\neq 0}\tf_\chi^{N_\chi}\right)\\
&\geq& -\frac{1}{\ln (k+N_c)}\left(\sum_{\chi\neq 0,1} N_\chi \ln(E_T (E_T\ovp)^{l(\chi,T)})\right)\\
&\geq& \frac{l(X,T)}{\ln (k+N_c)}\left(
- \frac{1}{l(X,T)}\left(\sum_{\chi\neq 0,1} N_\chi \ln E_T\right)
- \ln \frac{E_T l(X,T)\ln(k+N_c)}{N_c}
\right)\\
&\geq& \frac{l(X,T)}{1 + \ln N_c}\left(- \ln E_T + \ln N_c - \ln (E_T^2 k \ln(k+N_c))
\right)\\
&\geq& \frac{l(X,T)}{1 + \ln M^{1/\eps}}
\left(\ln M^{1/\eps} - 4\ln M - \ln(1+\ln M^{1/\eps})\right)\\
&\geq& l(X,T)(1 - 5\eps), 
\end{eqnarray*}
for $M$ large enough.
\end{proof}

\begin{proof} \textbf{(Proposition~\ref{prop:1})}
Let $T^*$ be a maximum likelihood tree with corresponding edge probabilities
$\bfp^*$, and $T^{**}$, be a maximum parsimony tree. Assume we have a 
polynomial-time algorithm which is guaranteed to return a tree $T'$
and edge probabilities $\bfp'$ such that
\begin{equation*}
-\ln \wtprob[X_0\,|\,T',\bfp']
\leq (1+c)\left(-\ln\wtprob[X_0\,|\,T^*,\bfp^*]\right).
\end{equation*}
Then the claims above and the optimality of $T^*$ imply that,
if $\bfp^{**}$ is chosen as in Claim~\ref{claim:1} (for $T=T^{**}$),
\begin{eqnarray*}
l(X,T')
&\leq& \frac{1}{1-5\eps}\left(\frac{-\ln \wtprob[X_0\,|\,T',\bfp']}{\ln(k+N_c)}\right)\\
&\leq& \frac{1+c}{1-5\eps}\left(\frac{-\ln \wtprob[X_0\,|\,T^*,\bfp^*]}{\ln(k+N_c)}\right)\\
&\leq& \frac{1+c}{1-5\eps}\left(\frac{-\ln \wtprob[X_0\,|\,T^{**},\bfp^{**}]}{\ln(k+N_c)}\right)\\
&\leq& \frac{(1+c)(1+2\eps)}{1-5\eps}l(X,T^{**})\\
&\leq& (1+c')l(X,T^{**}),
\end{eqnarray*}
for $\eps$ small enough.
\end{proof}

\begin{proof} \textbf{(Proposition~\ref{prop:2})}
Wareham~\cite[Theorem 45 Part 3]{Wa93} gives a reduction from 
vertex cover with bounded degree $B$ ($B$-VC) to maximum parsimony. (Wareham
actually defines MP as a Steiner tree problem on the Hamming cube $\{0,1\}^k$ but
the correspondence with our definition is straightforward.)
The reduction is such that
the existence of a $(1+c')$-approximation algorithm
for maximum parsimony implies the existence
of a $(1+2Bc')$-approximation algorithm
for $B$-VC. By~\cite{CT99}, for a sufficiently large $B$,
there is no $1.16$-approximation algorithm for $B$-VC unless
$\PP=\NP$. 
\end{proof}

\section*{Acknowledgment}

I thank Elchanan Mossel, Mike Steel and Tandy Warnow for discussions and encouragements. 
I gratefully acknowledge the partial support of 
CIPRES (NSF ITR grant \# NSF EF 03-31494), NSERC, NATEQ, and a Lo\`eve Fellowship. 
This project was initiated at a CIPRES retreat. I also thank Martin Nowak and the
Program for Evolutionary Dynamics at Harvard University where part of this work was done.


\begin{thebibliography}{1}

\bibitem{AC04}
L. Addario-Berry, B. Chor, M. T. Hallett, J. Lagergren, A. Panconesi, and T. Wareham, 
``Ancestral Maximum Likelihood of Evolutionary Trees is Hard,''
\emph{Journal of Bioinformatics and Computational Biology,} 
vol. 2, no. 2, pp. 257-271, 2004.

\bibitem{AB+99}
  	
R. Agarwala,
V. Bafna,
M. Farach,
B. Narayanan,
M. Paterson,
and M. Thorup,
``On the approximability of numerical taxonomy (fitting distances by tree metrics),''
SIAM Journal on Computing,
vol. 28, pp. 1073-1085, 1999.

\bibitem{AC99}
G. Ausiello, P. Crescenzi, G. Gambosi, V. Kann, A. Marchetti-Spaccamela, and M. Protasi, 
\emph{Complexity and Approximation,}
Springer, Berlin, 1999. 

\bibitem{Ca78}
J. Cavender,
``Taxonomy with confidence,'' 
\emph{Mathematical Biosciences,} vol. 40, pp. 271-280, 1978.

\bibitem{CT05}
B. Chor and T. Tuller,
``Maximum Likelihood of Evolutionary Trees is Hard,''
In: Proceedings of the 9th International Conference on 
Computational Molecular Biology (RECOMB 2005),
ACM Press, Cambridge, 2005.

\bibitem{CT99}
A. Clementi, and L. Trevisan,
``Improved Non-Approximability Results for Minimum Vertex Cover with Density Constraints,'' 
\emph{Theor. Comput. Sci.,} vol. 225, no. 1-2, pp. 113-128, 1999.

\bibitem{DJS86}
W. Day, D. Jonhson, and D. Sankoff,
``The computational complexity of inferring rooted phylogenies by parsimony,''
\emph{Mathematical Biosciences,}
vol. 81, pp. 33-42, 1986.

\bibitem{DS86}
W. Day and D. Sankoff,
``The computational complexity of inferring phylogenies by compatibility,''
\emph{Systematic Zoology,}
vol. 35, pp. 224-229, 1986.

\bibitem{EC64}
A. W. F. Edwards and L. L. Cavalli-Sforza,    
``Reconstruction of evolutionary trees,'' 
In: \emph{Phenetic and Phylogenetic Classification,} 
eds. V. H. Heywood and J. McNeill, 
Systematics Association, London, vol. 6, pp. 67-76, 1964.

\bibitem{Fa73}
J. S. Farris, 
``A probability model for inferring evolutionary trees,'' 
\emph{Systematic Zoology,} vol. 22,
pp. 250-256, 1973.

\bibitem{Fe81}
J. Felsenstein, 
``Evolutionary trees from DNA sequences: a maximum likelihood approach,''
\emph{J. Molecular Evolution,}  
vol. 17, pp. 368-376, 1981.
     
\bibitem{Fe03}
J. Felsenstein, 
\emph{Inferring Phylogenies,}
Sinauer Associates, Sunderland, 2004.

\bibitem{FG82}
L. Foulds and R. Graham,
``The Steiner problem in phylogeny is NP-complete,''
\emph{Advances in Applied Mathematics,}
vol. 3, pp. 43-49, 1982.

\bibitem{GJ76}
M. R. Garey and D. S. Johnson,  
\emph{Computers and Intractability. A Guide
to the Theory of NP-completeness,} 
W. H. Freeman, San Francisco, 1976.

\bibitem{Ne71}
J. Neyman,
``Molecular studies of evolution: a source of novel statistical problems,'' 
In: \emph{Statistical decision theory and related topics,} 
eds. S.S Gupta and J. Yackel, Academic Press, New York, pp. 1-27,  
1971. 

\bibitem{SS03}
C. Semple and M. Steel, 
\emph{Phylogenetics,} 
Oxford University Press, 2003.

\bibitem{TS97}
C. Tuffley and M. Steel, 
``Links between maximum likelihood and maximum parsimony under a simple model of site substitution,''
\emph{Bull Math Biol.,} vol. 59, no. 3, pp. 581-607,
1997.

\bibitem{Wa93}
H. T. Wareham,  
\emph{On the Computational Complexity of Inferring Evolutionary Trees,} 
M.Sc. thesis, Technical Report no. 9301, 
Department of Computer Science, Memorial University of Newfoundland,
1993.

\end{thebibliography}
\end{document}